\documentstyle[12pt,titlepage,twoside]{article}
\font\twelvegtc=eufm10 scaled 1200
\font\tengtc=eufm10
\font\ninegtc=eufm9
\font\sevengtc=eufm7
\font\fivegtc=eufm5
\newfam\gtcfam
\def\gtc{\fam\gtcfam\twelvegtc}
\textfont\gtcfam=\twelvegtc
\scriptfont\gtcfam=\ninegtc
\scriptscriptfont\gtcfam=\sevengtc
\font\twelveBBB=msbm10 scaled 1200
\font\tenBBB=msbm10
\font\sevenBBB=msbm7
\newfam\BBBfam
\def\BBB{\fam\BBBfam\twelveBBB}
\textfont\BBBfam=\twelveBBB
\scriptfont\BBBfam=\tenBBB
\scriptscriptfont\BBBfam=\sevenBBB
\def\QQ{{\BBB Q}}   \def\ZZ{{\BBB Z}}    \def\NN{{\BBB N}}
\def\RR{{\BBB R}}   \def\CC{{\BBB C}}    \def\PP{{\BBB P}}
\def\AA{{\BBB A}}   \def\aA{{\cal A}}    \def\kK{{\cal K}}
\def\pP{{\cal P}}   \def\mM{{\cal M}}    \def\oO{{\cal O}}
\def\VEC#1,#2{(#1_1,#1_2,\dots,#1_{#2})}
\def\OVEC#1,#2{(#1_0,#1_1,#1_2,\dots,#1_{#2})}
\def\SET#1,#2{\{#1_1,#1_2,\dots,#1_{#2}\}}
\def\OSET#1,#2{\{#1_0,#1_1,\dots,#1_{#2}\}}
\def\FAM#1,#2{\ #1_1,#1_2,\dots,#1_{#2}\ }
\def\BSER#1,#2,#3{#1_0+#1_1#2+\cdots+#1_{#3}#2^{#3}}
\def\SER#1,#2{#1_0+#1_1#2+#1_2#2^2+\cdots}
\def\POL#1,#2,#3{#1_0#2^{#3}+#1_1#2^{#3-1}+\cdots+#1_{#3-1}#2+#1_{#3}}
\def\UPOL#1,#2,#3{#2^{#3}+#1_1#2^{#3-1}+\cdots+#1_{#3-1}#2+#1_{#3}}
\def\Ddots{\mathinner{\mkern1mu\raise1pt\hbox{.}\mkern2mu\raise4pt\hbox
                  {.}\mkern2mu\raise7pt\vbox{\kern7pt\hbox{.}}\mkern1mu}}
\def\bydef{\stackrel{\rm def}{=}}
\def\im{{\rm im}\,}
\def\rk{{\rm rk}\,}
\def\ann{{\rm Ann}\,}
\def\tr{{\rm tr}\,}
\def\ch{{\rm ch}\,}
\def\Id{{\rm Id}}
\def\rp#1,#2{{{\rm Rp}_{#1}}#2}
\def\lp#1,#2{{{\rm Lp}_{#1}}#2}
\def\RM#1,#2{{{\rm Rm}_{#1}}#2}
\def\LM#1,#2{{{\rm Lm}_{#1}}#2}
\def\bl{\!\in\!}
\def\Is{{\sl Isom}\,}
\def\hom{{\rm Hom}}
\def\ext{{\rm Ext}}
\def\fa#1{\;\forall#1\;}
\def\fain#1,#2{\fa#1\bl#2}
\def\map#1,#2{#1\longrightarrow#2}
\def\MAP#1,#2,#3{#1\,\colon\;\map#2,#3}
\def\sp#1,#2{\langle#1,#2\rangle}
\def\cc#1,#2{\left[\frac{\textstyle #1}{\textstyle #2}\right]}
\def\openrow#1#2#3{\setbox0=\vbox{\hbox
    {\vrule height#2 width#3\kern#2\vrule height#2 width0pt}\hrule height#3}
    \hbox{\leaders\copy0\hskip#1\wd0\vrule width#3}}
\def\row#1#2#3{\vbox{\hrule height#3\openrow{#1}{#2}{#3}}}
\def\Yr#1{\row{#1}{1.5ex}{.1ex}}
\def\yr#1{\openrow{#1}{1.5ex}{.1ex}}
\def\DY#1\endDY{\baselineskip=1ex\lineskip=0pt\lineskiplimit=0pt{\vcenter
    {\Yr#1}}}
\def\openclm#1#2#3{\setbox0=\vbox{\hrule height#3\hbox
    {\vrule width0pt\kern#2\vrule width#3 height#2}}\vtop
    {\leaders\copy0\vskip#1\ht0\hrule height#3}}
\def\clm#1#2#3{\hbox{\vrule width#3\openclm{#1}{#2}{#3}}}
\def\Yc#1{\clm{#1}{1.5ex}{.1ex}}
\def\yc#1{\openclm{#1}{1.5ex}{.1ex}}
\def\CDY#1\endCDY{{\vcenter{\hbox{\Yc#1}}}}
\newcounter{No}
\newcounter{SubNo}[No]
\newcounter{SubSubNo}[SubNo]
\newcounter{cond}

\renewcommand{\thecond}{(\Alph{cond})}
\renewcommand{\theequation}{f.\arabic{No}.\arabic{SubNo}\Alph{equation}}
\renewcommand{\theNo}{\S\arabic{No}}
\renewcommand{\theSubNo}{\arabic{No}.\arabic{SubNo}}
\renewcommand{\theSubSubNo}{\arabic{No}.\arabic{SubNo}.\arabic{SubSubNo}}
\def\No#1{\refstepcounter{No}\par\vfill\pagebreak[3]\noindent
                               {\large\bf\theNo.\hspace{2pt}#1.}\par}
\def\SubNo#1{\refstepcounter{SubNo}\setcounter{equation}{0}\vspace{2ex}\par\noindent
                               {\bf\theSubNo.\hspace{2pt}#1.}\hspace{1ex}}
\def\SubSubNo#1{\refstepcounter{SubSubNo}\vspace{1ex}\par\noindent
                               {\bf\theSubSubNo.}\hspace{2pt}#1.\hspace{1ex}}
\def\EF{\endgroup\par\vspace{1ex}\par}
\def\EP{\par\nopagebreak[3]\noindent$\Box$\par}
\def\Pr{\SubSubNo{PROPOSITION}\begingroup\sl}
\def\Df{\SubSubNo{DEFINITION}\begingroup\sl}
\def\Cj{\SubSubNo{CONJECTURE}\begingroup\sl}
\def\Cl{\SubSubNo{COROLLARY}\begingroup\sl}
\def\proof{\par\noindent{\sc Proof.}\hspace{1ex}}
\def\Cs{\setcounter{cond}{0}\begin{description}}
\def\cond{\refstepcounter{cond}\item{\thecond\hspace{2ex}}}
\def\EC{\end{description}}
\title{Non-symmetric orthogonal geometry of Grothendieck rings of coherent
sheaves on projective spaces}
\author{A.L.Gorodentsev\thanks{This paper was started at the University of
Stockholm on June 1994 and finished at the Max-Plank-Institut f\"ur Mathematik
on August 1994}\thanks{In Moscow autor is supported by
the foundation PRO MATHEMATICA (France) and the J.Sorros foundation (USA)}\\
Algebra Section of the Steklov Mathematical Institute\\
GSP-1 Vavilova 42, Moscow, Russia \\
e-mail: gorod@alg.mian.su}
\date{June-August 1994}
\begin{document}
\maketitle
\begin{abstract}
In this paper we consider orthogonal geometry of the free $\protect\ZZ$-module
$K_0(\protect\PP_n)$ with respect to the non-symmetric unimodular bilinear form
$$\chi(E,F)=\sum (-1)^\nu\dim\ext^\nu(E,F).
$$
We calculate the isometry group of this form and describe invariants of its
natural action on $K_0(\protect\PP_n)$. Also we consider some general
constructions with non-symmetric unimodular forms. In particular, we discuss
orthogonal decomposition of such forms and the action of the braid group on a
set of semiorthonormal bases. We formulate a list of natural arithmetical
conjectures about semiorthogonal bases of the form $\chi$.
\end{abstract}
\markboth{A.L.Gorodentsev}{Non-symmetric Orthogonal Geometry of
$K_0(\protect\PP_n)$}
%
%
\No{Introduction}
\SubNo{The helix theory and the problem of description of exceptional sheaves
on $\PP_n$}
The helix theory is a cohomology technique to study derived categories of
coherent sheaves on some algebraic varieties. It appears first in [GoRu] and
[Go1] as the way to construct the {\it exceptional bundles\/} on $\PP_n$, i.e.
locally free sheaves $E$ such that
$$\dim \ext^0(E,E)=1,\;\ext^i(E,E)=0\;\forall\,i\!\ge\!1
$$
Since then the helix theory was developed in the context of general
triangulated categories in [Go2],[Go3],[Bo1],[Bo2].[BoKa]. The main idea of
this theory is to consider {\it exceptional bases\/} of a triangulated
category, i.e. collections of objects $\{E_0,E_1,\dots,E_n\}$ that generate the
category and have the following properties
$$\dim\hom^0(E_i,E_i)=1\,,\quad\hom^\nu(E_i,E_i)=0\;\forall\nu\!\not=\!1
$$
$$\hom^\mu(E_i,E_j)=0\;\forall\mu\;{\rm and}\;\forall i\!>\!j.
$$
The simplest example of a such collection is the collection
$$\{{\cal O}, {\cal O}(1),\dots,{\cal O}(n)\}
$$
of invertible sheaves on $\PP_n$. The main problem is to describe all such
collections. The most important fact in the study of this problem is that there
exists an action of the braid group on the set of exceptional collections of a
given length. Transformations of exceptional collections by generators of the
braid group are called {\it mutations\/}. The mutations make possible to
construct an infinite set of exceptional collection starting from a given one
(see [Go1], [Go2]).

Thus a description of all exceptional sheaves on $\PP_n$ splits into
three steps. We have to prove the following three conjectures.
\Cj
Any exceptional object in the bounded derived category of coherent sheaves on
$\PP_n$ is quasiisomorphic to a shifted image of an exceptional locally free
sheaf.
\EF
\Cj
Any exceptional collection (in particular, each exceptional object itself)
in the bounded derived category of coherent sheaves on $\PP_n$ can be included
in an exceptional basis of the derived category.
\EF
\Cj In bounded derived category of coherent sheaves on $\PP_n$ the braid group
acts transitively on exceptional collections of any given length.
\EF
\noindent
All these three conjectures hold on $\PP_2$ (see [GoRu], [Go2]), and the third
conjecture holds on $\PP_3$ for exceptional collections of maximal length
(generating the derived category, -- see [No]). Discussion of these problems in
the general context and the survey of corresponding results see in [Go4].
\SubNo{Subject of this paper}
In this paper we consider an arithmetical analog of problems formulated above.
Let us consider the Grothendieck group $K_0(\PP_n)$ as a free $\ZZ$-module of
finite rank with non-symmetric unimodular bilinear form
$$\chi(E,F)=\sum (-1)^\nu\dim\ext^\nu(E,F).
$$
\Df A collection of vectors $\SET e,k\subset K_0(\PP_n)$ is called {\it
exceptional\/} or {\it semiorthonormal\/} if the Gram matrix of the form $\chi$
at this collection is upper triangular with units on the main diagonal.
\EF
Obviously, any exceptional collection of sheaves produces an
se\-mi\-or\-tho\-nor\-mal collection of vectors in $K_0$. If we are going to
work only in terms of $K_0$ and $\chi$, then we can not distinguish such
collections and their images
with respect to the action of isometries of the form $\chi$.
\Df
$\ZZ$-linear operator $\varphi\,\colon\;K_0(\PP_n)\longrightarrow K_0(\PP_n)$
is called {\it isometric\/} if $\chi(v,w)=\chi(\varphi v,\varphi w)\;\forall
v,w$. The group of all isometric operators is denoted by $\Is$ and is called
the {\it isometry group\/}.
\EF
In \ref{formclas},\ref{Pn} we will prove that this group is an unipotent
Abelian algebraic group of dimension $[(n+1)/2]$. It has two connected
components, and the component of the identity is a direct sum of standard
1-dimensional additive groups. We write  explicit formulas for the natural
action of isometries on $K_0(\PP_n)$ and describe invariants of this action.
All this may be considered as the first step in the direction of the following
conjecture.

\Cj
A vector $e$ such that $\chi(e,e)=1$ represents (up to the action of
isometries) a class of exceptional sheaf if and only if it can be included in
some semiorthonormal basis of $K_0$.
\EF

In \ref{nonsyg} we consider some general constructions of non-symmetric
orthogonal geometry. In particular, we define an action of the braid group on
the set of all semiorthonormal collections of a given length, and introduce the
notions of the {\it canonical operator\/} and the {\it canonical algebra\/} of
given bilinear form, which play an important role in the general classification
of bilinear non-degenerate forms. We discuss this classification (over an
algebraically closed field of characteristic zero) in \ref{formclas} and give a
general approach to the problems like ones considered in [Ru].

This part of paper is the first little step in direction of the following
conjecture, which reformulate the main conjecture of helix theory
in terms of linear algebra.

\Cj Any semiorthonormal basis of $K_0(\PP_n)$ may be obtained from any other
one by changing signs of basic vectors and the action of braid group and
isometries.
\EF
%
%
\No{Non-symmetric orthogonal geometry}\label{nonsyg}
\SubNo{Notations}
Let $M$ be a free \ZZ-module of a finite rank equipped with an integer bilinear
form $M\times M\longrightarrow\ZZ$, which we denote by
$$ \sp*,*\,\colon\;v,w\mapsto\sp v,w.
$$
Submodules of $M$ will be usually denoted by $U,V,W,\dots$ Restriction of the
bilinear form onto submodule $U\!\subset\!M$ is denoted by $\sp*,*_U$.

For a given submodule $U\!\subset\!M$ the submodules
$$U^\bot=\{w\!\in\!M | \sp w,u=0\;\forall u\!\in\!U\}
$$
$$^\bot U=\{w\!\in\!M | \sp u,w=0\;\forall u\!\in\!U\}
$$
are called {\it right\/} and {\it left\/} orthogonals of $U$. So, $\sp
U,{U^\bot}=\sp {^\bot U},U=0$.

Fixing some basis $e=\SET e,n\!\subset\!M$ we denote by $\chi$ or $\chi(e)$ the
corresponding Gram matrix (the element $\chi_{ij}=\sp {e_i},{e_j}$ is placed in
$i$-th row and $j$-th column of this matrix).

We denote by $M^*$ the dual module $\hom_\ZZ(M,\ZZ)$.
\SubNo{Unimodularity and correlations}
For any given bilinear form on $M$ we can consider two linear operators:

{\it left correlation\/} $\MAP\lambda,M,{M^*}\,\colon\;v\mapsto\sp v,*$

\noindent and

{\it right correlation\/}  $\MAP\varrho,M,{M^*}\,\colon\;v\mapsto\sp*,v$.

\noindent The bilinear form on $M$ is uniquely determinated by each of them. In
fact, the Gram matrix $\chi$ coincides with a matrix of right correlation
written with respect to a pair of dual bases of $M$ and $M^*$, and a matrix of
the left correlation is its transpose. Hence, we get the following proposition:
\Pr The conditions:
\Cs
\cond left correlation is an isomorphism;
\cond right correlation is an isomorphism;
\cond $\det\chi=\pm1$;
\EC
are pairwise equivalent.
\EF\EP

We call a bilinear form on $M$ to be {\it unimodular\/} if it satisfies the
above conditions.

{\it In this paper we will always assume  that the form on $M$ is
unimodular.\/}

Note that if we identify $M$ with $M^{**}$ in the usual way, then the dual
operator to each of two correlations coincides with the other correlation:
$$\MAP{\varrho^*},{M^{**}=M},{M^*}\quad\mbox{ is equal to }
                                                    \quad\MAP\lambda,M,{M^*}
$$
$$\MAP{\lambda^*},{M^{**}=M},{M^*}\quad\mbox{ is equal to }\quad
                                                         \MAP\varrho,M,{M^*}
$$
\SubNo{Canonical operator}
Using correlations, we associate with any unimodular bilinear form on $M$ a
linear operator
$$\MAP{\kappa=\varrho^{-1}\lambda},M,M.
$$
This operator is called {\it canonical\/}. It is uniquely determinated by the
following condition:
$$\sp v,w=\sp w,{\kappa v} \fain{v,w},M .
$$
A matrix of $\kappa$ in any basis $e$ of $M$ is expressed in terms of the Gram
matrix $\chi=\chi(e)$ by the formula
$$\kappa=\chi^{-1}\chi^t
$$
Note that the map $\chi\mapsto\chi^{-1}\chi^t$ is equivariant with respect to
the standard actions of linear automorphisms of $M$ on bilinear forms and on
linear operators.

The canonical operator is isometric:
$$\sp v,w=\sp w,{\kappa v}=\sp{\kappa v},{\kappa w} \fain{v,w},M
$$
\SubNo{Dual operators and canonical algebra}\label{ca}
One can associate with any linear operator $\MAP\varphi,M,M$ a pair of its
dual.
They are uniquely determinated by the following conditions:

the {\it  left dual operator\/} $^\vee\varphi$:
 $\sp{^\vee\varphi v},w=\sp v,{\varphi w}$

the {\it  right dual operator\/} $\varphi^\vee$:
$\sp v,{\varphi^\vee w}=\sp{\varphi v},w.$

\noindent Expressions for their matrices are given by
$$^\vee\varphi=(\chi^{-1})^t\varphi^t\chi^t\;;\qquad\varphi^\vee=\chi^{-1}
                                                                \varphi^t\chi.
$$
In general, $^\vee\varphi\ne\varphi^\vee$ for a non-symmetric form on $M$. But
direct computation shows that the following proposition holds.
\Pr
The conditions
\Cs
\cond $^\vee\varphi=\varphi^\vee$
\cond $^{\vee\vee}\varphi=\varphi$
\cond $\varphi^{\vee\vee}=\varphi$
\cond $\varphi\kappa=\kappa\varphi$
\EC
are pairwise equivalent.
\EF\EP
An operator $\varphi$ is called {\it reflexive\/} if
$\varphi^{\vee\vee}=\varphi$ .
Reflexive operators form a subalgebra in the algebra of all linear
endomorphisms of $M$. This algebra coincides with the centralizer of the
canonical operator. We call it a {\it canonical algebra\/} of a bilinear form
on $M$ and denote by
$\aA$.

Obviously, the following sets of operators belong to $\aA$:

 $\aA^+=\{\varphi\,|\;^\vee\varphi=\varphi=\varphi^\vee\}$ --- the submodule of
{\it selfdual\/} operators;

 $\aA^-=\{\varphi\,|\;^\vee\varphi=-\varphi=\varphi^\vee\}$ --- the submodule
of {\it antiselfdual\/} operators;

 $\Is=\{\varphi\,|\;^\vee\varphi=\varphi^{-1}=\varphi^\vee\}$ --- the subgroup
of {\it isometric\/} operators.

Let us now consider the vector space $M_\QQ=\QQ\otimes_\ZZ M$ and the
corresponding canonical algebra $\aA_\QQ=\QQ\otimes_\ZZ\aA$. Evidently,
$\aA_\QQ$
splits into direct sum of the subspaces of selfdual and antiselfdual operators
$$\aA_\QQ=\aA^+_\QQ\oplus\aA^-_\QQ.
$$
$\Is$ is an algebraic group in the sense of [Se] and the arguments of [Se]
(ch.1,th.5) give us immediately the following proposition.
\Pr Subspace of all antiselfdual operators coincides with the Lie algebra of
algebraic group of all isometric operators:
$${\cal L}ie (\Is_\QQ)=\aA^-_\QQ
$$
\EF\EP
On the other side, we have
\Pr\label{equviv} Two bilinear forms $\sp*,*_1$ and $\sp*,*_2$ on a given
\ZZ-module $M$ have the same canonical operator if and only if there exists an
operator $\psi$, which is selfdual with respect to both forms and satisfies the
identity
$$\sp v,w_1=\sp v,{\psi w}_2\fain{v,w},M.
$$
\EF
\proof Taking the dual to the identity
$$\kappa=\varrho_1^{-1}\lambda_1=\varrho_2^{-1}\lambda_2
$$
we get
$$\kappa^*=\varrho_1\lambda_1^{-1}=\varrho_2\lambda_2^{-1}.
$$
Hence, $\psi=\varrho_2\varrho_1^{-1}=\lambda_2\lambda_1^{-1}$ satisfies the
identities
$$\sp v,{\psi w}_2=\sp v,w_1=\sp{\psi v},w_2.
$$
\EP
\Cl There exists the 1-1 correspondence between unimodular bilinear forms,
which have same canonical operators, and unimodular operators, which are
selfdual with respect to any one of these forms.
\EF\EP
\SubNo{Admissible submodules, orthogonal projections and mutations}
All constructions of this and two following sections are trivial reformulations
of the technique of orthogonal decomposition in triangulated categories (see
[Go2] and [BoKa1]).

A submodule $U\!\subset\!M$ is called {\it admissible\/} if it satisfies any of
the following equivalent conditions:
\Cs
\cond the restricted form $\sp*,*_U$ is unimodular;
\cond there exists a linear projection $\MAP{\lp U,{}},M,U$ such that
$$\sp v,u_M=\sp{\lp U,v},u_U\fain u,U\fain v,V;
$$
\cond there exists a linear projection $\MAP{\rp U,{}},M,U$ such that
$$\sp u,v_M=\sp u,{\rp U,v}_U\fain u,U\fain v,V;
$$
\cond $M=U\oplus{^\bot U}$;
\cond $M=U^\bot\oplus U$;
\EC
\noindent (equivalences (A)$\Leftrightarrow$(B)$\Leftrightarrow$(D) and
 (A)$\Leftrightarrow$(C)$\Leftrightarrow$(E) are standard in linear algebra
and we omit proofs).

Operators $\lp U,{}$ from (B) and $\rp U,{}$ from (C) are called {\it left\/}
and {\it right orthogonal projections onto $U$\/} respectively. Note that they
are left and right adjoint operators to the inclusion $U\hookrightarrow M$.

Of course, if a submodule $U\!\subset\!M$ is admissible, then both its
orthogonals $^\bot U$ and $U^\bot$ are admissible too. If a vector $v\bl M$ is
written in the form
$$v=v_{U^\bot}+v_{U}\;,\qquad\mbox{where }v_{U^\bot}\bl U^\bot\,,\;v_{U}\bl U,
$$
then obviously
$$\rp U,{(v)}=v_U\;,\quad \lp {U^\bot},{(v)}=v_{U^\bot}.
$$
Hence, we have the identities
$$v=\lp {U^\bot},{(v)}+\rp U,{(v)}=\lp U,{(v)}+\rp{^\bot U},{(v)}\;,
$$
which give decompositions of any vector $v\bl M$ as an element of
$U^\bot\oplus U$ and as an element of $U\oplus{^\bot U}$ respectively. It we
rewrite these identities in the form
$$\lp {U^\bot},{(v)}=\underbrace{ \lp U,{(v)}-\rp U,{(v)} }_{\mbox {from $U$} }
                      +\underbrace{ \rp {^\bot U},{(v)} }_{\mbox {from $ {^\bot
U} $} }
$$
or in the form
$$\rp{^\bot U},{(v)}=\underbrace{\lp {U^\bot},{(v)}}_  {\mbox{from $U^\bot$}}
                                  +\underbrace{\rp U,{(v)}-\lp
U,{(v)}}_{\mbox{from $U$}}                   $$
then we see that
$$\rp {^\bot U},{}{\scriptstyle\circ}\lp {U^\bot},{}=\rp {^\bot
U},{}\qquad\mbox{and}\qquad
               \lp {U^\bot},{}{\scriptstyle\circ}\rp {^\bot U},{}=\lp
{U^\bot},{}\;.
$$
Hence, the restriction of $\lp{U^\bot},{}$ onto the submodule $^\bot U$ and the
restriction of $\rp{^\bot U},{}$ onto the submodule $U^\bot$ are linear
isomorphisms between these submodules inverse to each other.

We call them {\it left mutation of $^\bot U$\/} and {\it right mutation of\/}
$U^\bot$ with respect to $U$ and denote by $\MAP{\LM{U},{}},{^\bot U},{U^\bot}$
and $\MAP{\RM{U},{}},{U^\bot},{^\bot U}$. They are well defined by the
properties:
$$\LM{U},{}\;\colon\,\rp {^\bot U},v\mapsto\lp{U^\bot},v\;\fain v,M
$$
$$\RM{U},{}\;\colon\,\lp {U^\bot},v\mapsto\rp{^\bot U},v\;\fain v,M.
$$
Moreover, direct computation shows that they are isometries. In fact,
$\forall\,u_1^\bot,u_2^\bot\in{U^\bot}$, which are decomposed in $U\oplus{^\bot
U}$ as
$u_\nu^\bot=u_\nu+{^\bot u_\nu}\quad(\nu=1,2)$ we have
\begin{eqnarray*}
\sp u_1^\bot,{u_2^\bot}&=&\sp u_1+{^\bot u_1},{u_2^\bot}=
            \sp {^\bot u_1},{u_2^\bot}=\sp {^\bot u_1},{u_2+{^\bot u_2}}=\\
                       &=&\sp {^\bot u_1},{{^\bot u_2}}.
\end{eqnarray*}
So, we get

\Pr\label{im} Left and right mutations with respect to admissible submodule
$U\!\subset\!M$ are isometric isomorphisms between ${^\bot U}$ and $U^\bot$
inverse to each other.
\EF\EP

\SubNo{Semiorthogonal direct sums}\label{cap}Let $M_1$ and $M_2$ be two modules
equipped with unimodular bilinear forms $\sp*,*_1$ and $\sp*,*_2$. Suppose that
 direct sum $M=M_1\oplus M_2$ is equipped with bilinear form $\sp*,*_M=\sp*,*$
such that $M_1=M_2^\bot$  (i.e. $\sp M_2,{M_1}=0$) and the restrictions of
$\sp*,*$ onto $M_1$ and $M_2$ coincide with $\sp*,*_1$ and $\sp*,*_2$
respectively. In this case we call $M$ to be a {\it semiorthogonal sum\/} of
$M_1$, $M_2$.

The form on $M$ is automatically unimodular too and we have
$$\sp u_1+u_2,{v_1+v_2}=\sp u_1,{v_1}_1+\sp u_1,{v_2}_M+\sp u_2,{v_2}_2\;,
$$
where the second summand can be expressed in two ways:
$$\sp u_1,{v_2}_M=\sp{\lp2,{u_1}},{v_2}_2=\sp u_1,{\rp1,{v_2}}_1
$$
(here $\lp2,{}=\lp M_2,{}|_{_{M_1}}$ and $\rp1,{}=\rp M_1,{}|_{_{M_2}}$ are the
orthogonal projections).

Hence, the structure of semiorthoghonal sum is uniquely determinated by any one
of two adjoint to each other linear operators
$$\MAP{\lp2,{}},{M_1},{M_2}\qquad\mbox{or}\qquad\MAP{\rp1,{}},{M_2},{M_1}\;.
$$

\Pr\label{canmatr}
Let $M=M_1\oplus M_2$ be a semiorthogonal sum. The canonical operator
$\kappa=\kappa_M$ is represented in terms of this decomposition by the
following matrix
$$\left(
    \begin{array}{cc}

\kappa_1-\rp1,{}{\scriptstyle\circ}\kappa_2{\scriptstyle\circ}\lp2,{}&-\rp1,{}\kappa_2\\
      \kappa_2\lp2,{}&\kappa_2
    \end{array}
 \right)\;,
$$
where $\kappa_1$ and $\kappa_2$ are the canonical operators on $M_1$ and $M_2$
respectively.
\EF
\proof
Let
$$\kappa_M=\left(
    \begin{array}{cc}
      \kappa_{11}&\kappa_{12}\\
      \kappa_{21}&\kappa_{22}
    \end{array}
 \right)\;,
$$
where $\MAP\kappa_{\mu\nu},M_\nu,{M_\mu}$. For any $u=u_1+u_2$, $v=v_1+v_2$
consider the identity
$$\sp u,v=\sp v,{\kappa u}\;.
$$
We write its left side as
\begin{eqnarray*}
   \sp u,v&=&\sp u_1,{v_1}_1+\sp{\lp2,u_1},{v_2}_2+\sp u_2,{v_2}_2=\\
          &=&\sp v_1,{\kappa_1 u_1}+\sp v_2,{\kappa_2\lp2,{(u_1)}}+\sp
v_2,{\kappa_2 u_2}
\end{eqnarray*}
and its right side as
$$\sp v,{\kappa u}=\sp v_1,{\kappa_{11} u_1+\kappa_{12} u_2}+
                   \sp v_1,{\rp1,{(\kappa_{21} u_1+\kappa_{22} u_2)}}+
                   \sp v_2,{\kappa_{21} u_1+\kappa_{22} u_2}\;.
$$
Comparing the components, we obtain
$$\kappa_1u_1=\kappa_{11}u_1+\kappa_{12}u_2+\rp1,{\kappa_{21}u_1}+\rp1,{\kappa_{22}u_2}
$$
$$\kappa_{2}\lp2,{u_1}+\kappa_{2}{u_2}=\kappa_{21}{u_1}+\kappa_{22}{u_2}
$$
It follows from the second equation that $\kappa_{22}=\kappa_{2}$ and
$\kappa_{21}=\kappa_{2}{\scriptstyle\circ}\lp2,{}$. Hence the first equation
may be rewritten as
$$\left\{\begin{array}{rcl}

\kappa_1&=&\kappa_{11}+\rp1,{}{\scriptstyle\circ}\kappa_2{\scriptstyle\circ}\lp2,{}\\
                 0&=&\kappa_{12}+\rp1,{}{\scriptstyle\circ}\kappa_2
          \end{array}
  \right.\;.
$$
\EP

\Df
A pair $(U,V)$ of submodules $U\!\subset\!M$, $V\!\subset\!M$ is called {\it
semiorthogonal\/} if both are admissible and $\sp v,u=0\;\fain{u,v},M$.
\EF

Let $W=U\oplus V$ be the direct sum of two admissible submodules forming a
semiorthogonal pair $(U,V)$. Of course, $W$ is admissible too and
$W^\bot=U^\bot\cap V^\bot$, $^\bot W={^\bot U}\cap{^\bot V}$.

\Pr In above notations $\lp W^\bot,{}=\lp U^\bot,{}\circ\lp V^\bot,{} $.
\EF
\proof
For any $w^\bot\bl W^\bot=U^\bot\cap V^\bot$ we have $\fain m,M$:
$$\sp m,{w^\bot}=\sp {\lp V^\bot,m},{w^\bot}=\sp {\lp U^\bot,{\lp
V^\bot,m}},{w^\bot}.
$$
On the other side, $\im(\lp U^\bot,{}\circ\lp V^\bot,{})\subset W^\bot$ (this
follows from the assumption that $U\!\subset\!V^\bot$).
\EP

By the same way we conclude

\Pr $\rp ^\bot W,{}=\rp ^\bot V,{}\circ\rp ^\bot U,{} $.
\EF\EP

These propositions give

\Cl\label{cm} In above notations
$$\LM W,{}=\LM U,{}\circ\LM V,{}
$$
$$\RM W,{}=\RM V,{}\circ\RM U,{}.
$$
\EF\EP

\SubNo{Semiorthonormal collections and Braid Group action}\label{mutations}
The simplest example of an admissible submodule is a $1$-dimensional free
submodule $\ZZ e\bl M$ generated by a vector $e$ with $\sp e,e=1$. In this case
we have
$$\begin{array}{rclrcl}
      \rp e,v&=&\sp e,v e\;,&\lp e,v&=&\sp v,e e\;,\\
      \LM e,v&=&\lp e^\bot,v=v-\sp e,v e\;,&\RM e,v&=&\rp ^\bot e,v=v-\sp v,e
e\;.
  \end{array}
$$

If $W={^\bot e}$, then it follows from (\ref{cap}) that the scalar product on
the semiorthogonal sum $M=\ZZ e\oplus W$ is uniquely determinated by the scalar
product $\sp*,*_W$ on $W$ and the vector $\ell=\lp W,e\in W$.

By (\ref{canmatr}), the canonical operator $\kappa=\kappa_M$ is represented in
terms of the decomposition $M=\ZZ e\oplus W$ by the matrix
$$\kappa_M=\left(
    \begin{array}{cc}
      \lambda&\psi\\
      \kappa_W(\ell)&\kappa_{W}
    \end{array}
 \right)\;,
$$
where
$$\MAP\psi,W,{\ZZ e}\,\colon\;w\mapsto-\sp e,{\kappa_2w}\cdot
e=\sp\ell,{\kappa_2w}_W\cdot e
$$
and
$$\lambda=1-\sp
e,{\kappa_2\ell}_M=1-\sp\ell,{\kappa_2\ell}_W=1-\sp\ell,\ell_W\;.
$$
On the other hand $\sp\kappa_M e,e=\sp\lambda e,e+\sp\kappa_W(\ell),e=\lambda$.
We get

\Pr\label{trkap}
$\tr(\kappa_M)=\tr(\kappa_W)+\sp\kappa_M e,e=\tr(\kappa_W)+1-\sp\ell,\ell$ and
$\sp\kappa_M e,e=1-\sp\ell,\ell$.
\EF\EP\pagebreak[3]

Of course, not only one vector but any semiorthonormal collection of $(k+1)$
vectors $\OSET e,k$ (i.e. such that
$\sp e_\mu,{e_\nu}=0\;\forall \mu>\nu\,,\;\sp e_\nu,{e_\nu}=1\:\forall\nu$
) generate an admissible submodule $W$. For a such submodule $W$ we have
$$\LM W,{}=\LM e_0,{}\circ\LM e_1,{}\circ\cdots\circ\LM e_k,{},
$$
$$\RM W,{}=\RM e_k,{}\circ\RM e_{k-1},{}\circ\cdots\circ\RM e_0,{}
$$

Let us consider a module $M$ generated by a semiorthonormal pair $(a,b)$.
Define the {\it left\/} and {\it right mutations\/} of this pair by the
formulas
$$L(a,b)\bydef(\LM a,b,a)=(b-\sp a,b a,a)
$$
$$R(a,b)\bydef(b,\RM b,a)=(b,a-\sp a,b b)
$$
Note that these mutations may be considered as results of two Gram-Schmidt
orthogonalisations applied to the non-semiorthogonal pair $(b,a)$ in two
possible ways.

It follows from (\ref{im}) that
$$RL(a,b)=LR(a,b)=(a,b).
$$

Now, let us consider a module $M$ generated by a semiorthonormal triple
$(a,b,c)$ and its submodule $W$ generated by $(a,b)$. In this case $^\bot W$ is
generated by $c$. We can calculate the mutation $\LM W,c$ in two ways: taking
in (\ref{cm}) $U$ generated by $a$, $V$ generated by $b$, or taking  $U$
generated by $\LM a,b$, $V$ generated by $a$. The results must be the same, and
we get some kind of the {\it triangle equation\/}:
$$\LM a,{\LM b,c}=\LM {\LM a,b},{\LM a,c}.
$$

To present last two identities in more conceptual form let us consieder a
submodule $W\!\subset\!M$ generated by semiorthonormal collection $\OSET e,k$.
Denote by $L_\nu$ and $R_\nu$ the operations, which change the pair
$(e_{\nu-1},e_\nu)$ by its left and right mutations respectively
($\nu=1,2,\dots,k$). From our identities it follows immediately the following
\Pr Operations $L_\nu$ and $R_\nu$ satisfy the identities:
$$\begin{array}{rclr}
   R_\nu L_\nu&=&L_\nu R_\nu\;=\;{\rm Id}\\
   L_\nu L_{\nu-1}L_\nu&=&L_{\nu-1}L_\nu L_{\nu-1}&{\rm
                                           for}\quad\,\nu=2,3,\dots,k \\
   L_\mu L_\nu&=&L_\nu L_\mu&{\rm for}\quad\mu,\nu\,\colon\;|\mu-\nu|>1
\end{array}
$$
\EF\EP
\Cl The braid group acts by left mutations of neighboring pairs on the set of
semiorthonormal bases of an admissible submodule.
\EF\EP
%
%
%
\No{Decomposition of bilinear forms via canonical operator}\label{formclas}
\SubNo{Notations}
In this paragraph we consider a vector space $V$ over an algebraically closed
field $K$ of characteristic $0$ equipped with a non-degenerate bilinear form
$\sp*,*$.

The pair $\{V,\sp*,*\}$ is called {\it decomposable\/}, if $V=U\oplus W$, the
restrictions $\sp*,*_U$ and $\sp*,*_W$ are non-degenerate (in particular,
$U\ne0$, $W\ne0$), and $U$, $W$ are {\it biorthogonal\/} to each other, i.e.
$\sp U,W=\sp W,U=0$.

We are going to classify in all indecomposable pairs  $\{V,\sp*,*\}$ up to
isometric isomorphisms. The following proposition shows that the answer may be
given in terms of the canonical operator of the form on $V$.

\Pr Let $\{V,\sp*,*_V\}$,  $\{U,\sp*,*_U\}$ be two spaces  with
non\--de\-ge\-ne\-ra\-te bilinear forms. They are isometrically isomorphic to
each other if and only if there exists an isomorphism $\MAP\psi,U,V$ such that
$\psi\kappa_V=\kappa_U\psi$, where $\kappa_U$ and $\kappa_V$ are the canonical
operators of the forms on $U$ and $V$.
\EF
\proof
We may assume that two different forms $\sp*,*_1$ and $\sp*,*_2$ on the same
vector space $V$ are given and that these two forms have the same canonical
operator $\kappa$. It is sufficient to prove that in this case there exists a
linear isomorphism $\MAP\varphi,V,V$ such that
$$\sp v,w_1=\sp\varphi v,{\varphi w}_2\;\fain {v,w},V.
$$
In (\ref{equviv}) we have seen that there exists a selfdual operator
$\MAP\psi,V,V$ such that
$$\sp v,w_1=\sp v,{\psi w}_2\;\fain {v,w},V.
$$
Over an  algebraically closed field $K$ of characteristic $0$ we can find a
polynomial $F(t)\bl K[t]$ such that the operator $\varphi\bydef F(\psi)$
satisfy the equation $\varphi^2=\psi$. Since $\varphi$ is selfdual too, we
obtain $\sp v,w_1=\sp\varphi v,{\varphi w}_2\;\fain {v,w},V$.
\EP

So, non-degenerate non-symmetric bilinear form over algebraically closed field
of characteristic zero is uniquely determinated by Jordan normal form ot its
canonical operator. We will describe the correspondence between the root
decomposition of the canonical operator and biorthogonal decomposition of
original bilinear form. These results are not new and actually they may be
extracted from classical books [HoPe] (Book 2, Ch.IX) and [Ma].

For any linear operator $\MAP\varphi,V,V$ we will usually denote by $\lambda$,
$\mu$, $\dots$ its eigenvalues and by  $\varphi_\lambda$, $\varphi_\mu$,
$\dots$ --- corresponding differerences  $\varphi-\lambda E$, $\varphi -\mu E$,
$\dots$  The root subspaces, which corresponds to these eigenvalues, will be
denoted by
$V_\lambda$, $V_\mu$, $\dots$ . So,
$$V_\mu=\bigcup_{n\in\NN}\ker\varphi^n_\mu.
$$

\SubNo{Decomposition of an isometry}
Let $\MAP\varphi,V,V$ be an isometric operator with eigenvalues $\lambda$,
$\mu$ and let
$$v_m\stackrel{\varphi_\lambda}{\mapsto}
       v_{m-1}\stackrel{\varphi_\lambda}{\mapsto}
        \cdots\stackrel{\varphi_\lambda}{\mapsto}
           v_0\stackrel{\varphi_\lambda}{\mapsto}v_{-1}=0
$$
$$w_k\stackrel{\varphi_\mu}{\mapsto}
       w_{k-1}\stackrel{\varphi_\mu}{\mapsto}
        \cdots\stackrel{\varphi_\mu}{\mapsto}
           w_0\stackrel{\varphi_\mu}{\mapsto}w_{-1}=0
$$
be any two Jordan chains for operators $\varphi_\lambda=\varphi-\lambda E$ and
$\varphi_\mu=\varphi -\mu E$. Then for any $0\le i\le m$, $0\le j\le k$ we have
$$\sp v_i,{w_j}=\sp \varphi v_i,{\varphi w_j}=
             \sp\lambda v_i+v_{i-1},{\mu w_j+w_{j-1}}
$$
and hence
$$(1-\lambda\mu)\sp v_i,{w_j}=\lambda\sp v_i,{w_{j-1}}+\mu\sp v_{i-1},{w_j}+
                 \sp v_{i-1},{w_{j-1}}.
$$
Using decreasing induction we obtain that for  $\lambda\mu\ne1$
$$\sp v_i,{w_j}=\sp w_j,{v_i}=0.
$$
We have proved
\Pr If $\lambda\mu\ne1$, then two root subspaces $V_\lambda$ and $V_\mu$ of any
isometry $\MAP\varphi,V,V$ are biorthogonal to each other.
\EF\EP
\Cl  Let $\MAP\varphi,V,V$ be an an isometry of a space $V$ equipped with
non-degenerate bilinear form. Then $V$ splits into biorthogonal direct sum of
subspaces $W_\mu$, where:

-- for $\mu=\pm1$ $W_\mu$ coincides with the root subspace $V_\mu$ of $\varphi$
and restriction of the original form onto $W_\mu$ is nondegenerate;

-- for $\mu\ne\pm1$ $W_\mu$ coincides with the direct sum of root subspaces
$V_\mu\oplus V_{\mu^{-1}}$ and the original form restricts trivially onto each
of these two root subspaces and induces a nondegenerate pairing between them.
\EF\EP

In order to clarify the action of $\varphi$ on the subspaces $W_\mu=V_\mu\oplus
V_{\mu^{-1}}$ we denote $V_\mu$ by $V_+$ and $V_{\mu^{-1}}$ by $V_-$. Let us
identify $V_-^*$ with $V_+$ using non-degenerate pairing $\sp V_+,{V_-}$. So,
we can consider the dual to a linear operator $\MAP f,V_+,V_+$ as the operator
$\MAP f^*,V_-,V_-$ defined by the formula
$$\sp v_+,{f^*v_-}=\sp fv_+,{v_-}\quad\fain v_+,{V_+}\;\fain v_-,{V_-}.
$$
Finally, consider two nilpotent operators
$$\MAP\varphi_+\bydef(\varphi -\mu E)|_{_{V_+}},V_+,{V_+}\;;
$$
$$\MAP\varphi_-\bydef(\varphi -\mu^{-1} E)|_{_{V_-}},V_-,{V_-}\;.
$$
In the case $\mu=\pm1$ we put $V_+=V_-=V_\mu$ and
$\varphi_+=\varphi_-=\varphi_\mu$.

\Pr\label{ortcond} The following formulae hold:
\Cs
\cond $\varphi_+^*=-\mu\varphi^{-1}|_{_{V_-}}\varphi_-$;
\cond $\ker(\varphi_-^k)=\ann(\im(\varphi_+^k))$;
\cond $\im(\varphi_-^k)=\ann(\ker(\varphi_+^k))$;
\EC\EF
\proof
The first formula is checked by direct computation.

The other two follow from the first one by using the fact that for any pair of
dual linear operators
 $\MAP f,V,V$ and $\MAP f^*,V^*,{V^*}$
we have $\ker(f)=\ann(\im(f^*))$ and $\im(f)=\ann(\ker(f^*))$. In our case one
have to put $f=\varphi_+^k$, $f^*=(-\mu)^k\varphi^{-k}|_{_{V_-}}\varphi_-^k$
and to note that $\ker(f^*)=\ker(\varphi_-^k)$ and
$\im(f^*)=\im(\varphi_-^k)$, because $(-\mu)^k\varphi^{-k}|_{_{V_-}}$ is an
isomorphism, which commute with $\varphi_-$.
\EP
\Cl\label{jord} Nilpotent operators $\varphi_+$ and $\varphi_-$ have the same
Jordan normal form (the same cycle type).
\EF
\proof
It is convenient to represent an Jordan basis of a given nilpotent operator $f$
on a space $W$ by the Young diagram like the following one:
$$\DY7\yr5\yr5\yr2\yr1\endDY\qquad\stackrel{f}{\leftarrow}
$$
The cells of this diagram are in 1-1 correspondence with the basic vectors of
Jordan basis and $f$ takes each cell to its left neighboring and takes the
cells from the first left column to zero.

In terms of such representation the sum $S_k(f)$ of lengths of the first left
$k$ columns is equal to the $\dim\ker(f^k)$. From the other hand, the number of
cells forming these $k$ columns coincides with the number of cells  forming a
basis of a direct complement to the subspace $\im(f^k)$ (these cells are at the
right side of the Young diagram). So,
$$S_k(f)=\dim\ker(f^k)=\dim W-\dim\im(f^k).
$$
It follows from above proposition that in our case
$$\dim\ker(\varphi^k_+)=\dim(V_-)-\dim\im(\varphi_-^k).
$$
Hence, $S_k(\varphi_+)=S_k(\varphi_-)\;\forall k$ and our operators have the
same Young diagram.
\EP
\SubNo{Decomposition of the canonical operator}
Suppose now that $V$ is indecomposable and apply the previous results to the
canonical operator $\varphi=\kappa$. We see that there are two cases.

In the first case, which we call for a moment {\it the $\mu$-case\/},
$$V=V_+\oplus V_-\,\quad\kappa|_{_{V_+}}=\mu
E+\eta_+\,\quad\kappa|_{_{V_-}}=\mu^{-1} E+\eta_-,
$$
where $\mu\ne\pm1$ and $\eta_+$, $\eta_-$ are nilpotent operators of the same
cycle type. The restrictions $\sp*,*_{V_+}$, $\sp*,*_{V_-}$ are zeros and
pairing $\sp V_+,{V_-}$ is non-degenerate.

In the second case, which we call for a moment {\it the $\varepsilon$-case\/}
$$\kappa=\varepsilon E+\eta,
$$
where $\varepsilon=\pm1$ and $\eta$ is nilpotent.

The exact description of Jordan normal form of $\eta_\pm$ and $\eta$ in these
two cases will be given in two consequent propositions below.

\Pr\label{mucase}
If $V$ is indecomposable and the $\mu$-case takes place, then $\eta_+$ and
$\eta_-$ have only one Jordan cycle of the same length, i.e. the Young diagrams
of $\eta_\pm$ are of the form $\DY7\endDY\;$.
\EF
\proof
Denote by $K_\pm$ the kernel subspaces $\ker(\eta_\pm)\!\subset\!V_\pm$. We fix
some Jordan basis for $\eta_+$ in $V_+$ and denote by $C_+$ the direct
complement to the image subspace $\im(\eta_+)\!\subset\!V_+$  induced by this
choice, so $V_+=C_+\oplus\im(\eta_+)$. Let
$$e^+_k\stackrel{\eta_+}{\mapsto}
       e^+_{k-1}\stackrel{\eta_+}{\mapsto}
        \cdots\stackrel{\eta_+}{\mapsto}
           e^+_0\stackrel{\eta_+}{\mapsto}e^+_{-1}=0
$$
be the Jordan chain of maximal length for $\eta_+$ (corresponding to the upper
row of the Young diagram of $\eta_+$), $L_+$ be its linear span, and $W_+$ be
the linear span of all others basic vectors, i.e. $V_+=L_+\oplus W_+$.

It follows from proposition \ref{ortcond} that the pairing $\sp C_+,{K_-}$ is
non-degenerate. We fix the basis of $K_-$, which is dual to the basis of $C_+$
fixed above, and consider the vector $e^-_0$ of this basis such that $\sp
e^+_k,{e^-_0}=1$ and $\sp e^+,{e^-_0}=0$ for all other basic vectors $e^+\bl
C_+$. It follows from (\ref{ortcond}) and (\ref{jord}) that automatically
$e^-_0\bl\im(\eta_-^k)$. Hence, this vector can be included in a Jordan cain
$$e^-_k\stackrel{\eta_-}{\mapsto}
       e^-_{k-1}\stackrel{\eta_-}{\mapsto}
        \cdots\stackrel{\eta_-}{\mapsto}
           e^-_0\stackrel{\eta_-}{\mapsto}e^-_{-1}=0.
$$

Each vector $e^-_j$ of this chain is determinated by $e^-_{j-1}$ not uniquely
but modulo $K_-$.
Since the pairing $\sp C_+,{K_-}$ is non-degenerate, we can modify this chain
(in the unique way!) in order to have $e^-_j\bl C_+^\bot\;\forall j\!\ge\!1$.
Denote by $L_-$ the linear span of the chain chosed in the such way. We have,
in particular,
$$\sp C_+\!\cap\!W_+,{L_-}=0.
$$

It is easy to check that $L_-$ is biorthogonal to $W_+$. Actually, any $w^+\bl
W_+$ can be written as $w^+=\eta_+^mc^+$, where $c^+\bl C_+\cap W_+$. Hence,
$\forall j$ we have:
$$\sp w^+,{e_j^-}=\sp \eta_+^mc^+,{e_j^-}=(-\mu)^m\sp
c^+,{\eta_-^m\kappa^{-m}e_j^-} =0,
$$
because $\sp c^+,{L_-}=0$ and $L_-$ is invariant under the action of $\kappa$
and $\eta_-$. Orthogonality in the opposite direction follows immediately:
$$\sp L_-,{w^+}=\sp w^+,{\kappa L_-}=\sp w^+,{ L_-}=0.
$$

Starting from the others basic vectors of $K_-$ we can construct in the same
way as above a direct decomposition $V_-=L_-\oplus W_-$ such that $\sp
L_+,{W_-}=\sp W_-,{L_+}=0$. Hence, the subspace $L_+\oplus L_-\subset V$ is a
biorthogonal direct summand in $V$. Since $V$ is indecomposable, we have
$V=L_+\oplus L_-$.
\EP
\Pr\label{epscase}
If $V$ is indecomposable and the $\varepsilon$-case takes place, then either
$\eta$ has the Young diagram of the form
$$\underbrace{\DY5\endDY\cdots\DY5\endDY}_{n+1}
$$
and $\varepsilon=(-1)^n$ or $\eta$ has the Young diagram of the form
$$\underbrace{\DY5\yr5\endDY\cdots\DY5\yr5\endDY}_{n+1}
$$
and $\varepsilon=(-1)^{n+1}$.
\EF
\proof
Let $\eta^n\ne0$ but $\eta^{n+1}=0$. It follows from proposition \ref{ortcond}
that the bilinear form
$$(v,w)\bydef\sp v,{\eta^nw}
$$
is well defined and non-degenerate on the factor-space $V/\ker(\eta^n)$. The
calculation:
\begin{eqnarray*}
(w,v)&=&\sp w,{\eta^nv}=\sp{\eta^nv},w=\sp v,{\left(\eta^\vee\right)^n\kappa
w}=
                                     (-\varepsilon)^n\sp
v,{\kappa^{-n}\eta^n\kappa w}=\\
     &=&(-\varepsilon)^n\sp v,{(\varepsilon E+\eta)^{1-n}\eta^n w}=
                                     (-\varepsilon)^n\varepsilon^{1-n}\sp
v,{\eta^n w}=\\
     &=&(-1)^n\varepsilon(v,w)
\end{eqnarray*}
shows that this form is symmetric for $\varepsilon=(-1)^n$ and is
skew-symmetric for  $\varepsilon=(-1)^{n+1}$.

In the first case we can find an orthonormal basis of $V/\ker(\eta^n)$ with
respect to this symmetric form. If we construct a Jordan basis for $\eta$ in
$V$ starting with this orthonormal basis of $V/\ker(\eta^n)$, then we will be
in a position to apply the arguments from the proof of the previous
proposition. Exactly as above we can modify the Jordan chain of maximal length
in such a way that its linear span will be detached as biorthogonal direct
summand. Hence, in this case $\kappa$ has only one Jordan cycle and its length
is modulo 2 different from the eigenvalue of $\kappa$.

In the second case we can decompose $V/\ker(\eta^n)$ with respect to symplectic
form $(*,*)$ into direct sum of standard 2-dimensional symplectic planes, which
are orthogonal to each other. If we fix a symplectic basis in one of these
planes, then, as above, we can construct Jordan chains ended in these two
vectors in such a way that its linear span will be detached as biorthogonal
direct summand. Hence, in this case $\kappa$ has only two Jordan cycle of the
same length  and this length is modulo 2 equal to the eigenvalue of $\kappa$.
\EP

\SubNo{Forms of type 2}
Note that the second case of proposition \ref{epscase} can be considered as a
particular case of the proposition \ref{mucase}. Namely, in both cases the
structure of the space $V$ is described by the following definition.

\Df Let $K$ be an arbitrary field, $V$ be a vector space of even dimension
$\dim V=2k$  over $K$, and  $\mu\bl K$, $\mu\ne(-1)^{k+1}$. Non-degenerate
indecomposable bilinear form on the space $V$ is called to be of the {\it type
2\/}, if
$$V=V_+\oplus V_-\;.\quad\kappa|_{_{V_\pm}}=\mu^{\pm1}E+\eta_\pm\;,
$$
where $\dim V_+=\dim V_-=k$, $\eta_\pm^k=0$, $\eta^{k-1}_\pm\ne0$, the pairing
$\sp V_+,{V_-}$ is non-degenerate, and both restrictions $\sp*,*_{V_+}$,
$\sp*,*_{V_+}$ are identically equal to zero.
\EF

It is easy to see that non-degenerate indecomposable forms of the type 2
actually exist for any $k\bl\NN$ and $\mu\bl K$. Moreover, from the proofs
given above one can extract some standard form for their Gram matrix.

\Cl Each non-degenerate indecomposable form of the type 2 has at some
appropriate basis the Gram matrix
$$\left(
  \begin{array}{ccccc|ccccc}
    &&&&&&&&&\mu\\
    &&&&&&0&&\mu&1\\
    &&0&&&&&\Ddots&\Ddots&\\
    &&&&&&\mu&1&&0\\
    &&&&&\mu&1&&&\\
   \hline
    &&&&1\,&&&&&\\
    &0&&1\,&&&&&&\\
    &&\Ddots&&&&&0&&\\
    &\,1&&0&&&&&&\\
    \,1&&&&&&&&&\\
  \end{array}
   \right)
$$
\EF
\proof
First consider a form of the type 2. In this case $V=V_+\oplus V_-$, $\dim
V_+=\dim V_-=k$. Let us fix an arbitrary Jordan basis
$$e^+_k\stackrel{\eta_+}{\mapsto}
       e^+_{k-1}\stackrel{\eta_+}{\mapsto}
        \cdots\stackrel{\eta_+}{\mapsto}
           e^+_0\stackrel{\eta_+}{\mapsto}e^+_{-1}=0
$$
for $\eta_+$ in $V_+$. Since the pairing $\sp
\ker(\eta_-^j),{V_+/\im(\eta_+^j)}$ is non-degenerate $\forall j$, we can find
a sequence of vectors $v_j\bl\ker(\eta_-^j)$ (where $j=0,\,1,\,\dots\,,\,k$)
such that
\begin{eqnarray*}
          \sp v_j,{e_{k-j}^+}&=&1\\
       \sp v_j,{e_{k-\nu}^+}&=&0\quad\mbox{for $\nu=0,\,1,\,\dots\,,\,(j-1)$.}
\end{eqnarray*}
Note that by (\ref{ortcond}) we have also
\begin{eqnarray*}
     \qquad\;\sp v_j,{e_{k-\nu}^+}&=&0\quad\mbox{for
$\nu=(j+1),\,(j+2),\,\dots\,,\,k$.}
\end{eqnarray*}
These relations determinate the scalar products in the opposite order too:
\begin{eqnarray*}
      \sp e_{k-\nu}^+,{v_j}&=&\sp v_j,{\kappa e_{k-\nu}^+}=
                                     \mu\sp v_j,{e_{k-\nu}^+}+\sp
v_j,{e_{k-\nu-1}^+}=\\
                           &=&\left\{\begin{array}{cl}
                                   \mu&\mbox{for $\nu=j$}\\
                                     1&\mbox{for $\nu=j-1$}\\
                                     0&\mbox{for all other $\nu$}\quad.
                                \end{array}\right.
\end{eqnarray*}
We see that the Gram matrix of the basis
$\{e_0^+,\dots,e_k^+,v_0^+,\dots,v_k^+\}$ has the form what is needed.

On the other hand, simple direct computation shows that the canonical operator
of this Gram matrix actually has two Jordan cycles of length $k$ with
eigenvalues $\mu$ and $\mu^{-1}$.
\EP

\SubNo{Forms of type 1}\label{typeone}
Actually, much more interesting for us are the forms, which satisfy the
remainder first condition from proposition \ref{epscase}.

\Df Non-degenerate indecomposable bilinear form on a space $V$ of dimension
$\dim V=n+1$ is called to be of the {\it type 1\/}, if
$$\kappa=(-1)^n E+\eta\;.
$$
where $\eta^{n+1}=0$, $\eta^n\ne0$.
\EF

We will show in \ref{Pn} that the form on $K_0(\PP_n)$ is of type 1. Let us
consider the forms of type 1 in a more details. We put in this section $K=\QQ$,
because this case will be used in the next paragraph, but actually all results
are true for any fild of characteristic zero.

We fix a vector space $V$ of dimension $\dim V=n+1$ over \QQ\ and denote the
sign $(-1)^n$ by $\varepsilon$. Since the canonical operator of a form of type
1 has the form
$\kappa=\varepsilon E +\eta$ and nilpotent operator $\eta$ has a Jordan chain
of length $n+1$, the centralizer of the canonical operator in ${\rm
End}_\QQ(V)$ (i.e. the {\it canonical algebra\/} $\aA$, see \ref{ca}) coincides
with the commutative subring
$$\aA=\QQ[\eta]/\eta^{n+1}\subset{\rm End}_\QQ(V)\;.
$$

In order to study the involution $\scriptstyle\vee$ (see \ref{ca}) it is more
convenient to choose an other generator of the canonical algebra. Namely, let
$$\zeta=\frac{\varepsilon\kappa-E}{\varepsilon\kappa+E}=
          \frac{1}{2}\varepsilon\eta(E+\frac{1}{2}\varepsilon\eta)^{-1}\;.
$$
We have
$$\kappa=\varepsilon\frac{1+\zeta}{1-\zeta}\;;\qquad
                           \eta=2\varepsilon(\zeta+\zeta^2+\cdots+\zeta^n)
$$
and $\kappa$, $\eta$, $\zeta$ are uniquely determinated by each other.
Obviously, $\zeta$ is nilpotent, has the same Jordan normal form as $\eta$,
$\ker(\eta^i)=\ker(\zeta^i)\,\forall i$, and $\QQ[\zeta]=\QQ[\eta]$. By
proposition \ref{ortcond} the involution $\scriptstyle\vee$ is uniquely
determinated by the condition $\eta^\vee=-\varepsilon\kappa^{-1}\eta$. Hence,
this involution coincides with the involution of the ring $\QQ[\zeta]$, which
takes $\zeta$ to $-\zeta$, and acts on $\aA=\QQ[\zeta]/\zeta^{n+1}$ by the rule
$$f(\zeta)^\vee=f(-\zeta)\quad\fain f,\QQ[\zeta]\;.
$$
In particular, the Gram matrix of any Jordan basis $\{e_i\}$ for $\zeta$ is
uniquely determinated by its right column by the simple rule:
\begin{equation}\label{jogr}
  \sp e_i,{e_j}=\sp \zeta^{n-i}e_n,{\zeta^{n-j}e_n}=
               \left\{\begin{array}{cl}
                         (-1)^j\varepsilon\sp e_{i+j-n},{e_n}&
                                                  \mbox{for}n\le(i+j)\le2n\\
                         0&\mbox{for}(i+j)<n
                       \end{array}
               \right.\;.
\end{equation}

We see that the subspace $\aA_+\!\subset\!\aA$ of all selfdual operators
coincides with the subspace of all operators $f(\zeta)$ represented by even
polynomials $f$. The subspace $\aA_-={\cal L}{\it ie}(\Is)\!\subset\!\aA$ of
all antiselfdual operators is generated by $k=\cc n+1,2$ odd powers
$\zeta,\,\zeta^3,\,\dots\,,\,\zeta^{2k-1}$.

\Pr The isometry group of a form of type 1 is Abelian and has two connected
components. The component of the identity is isomorphic to the direct product
of $\cc n+1,2$ standard 1-dimensional additive unipotent algebraic groups.
\EF
\proof
It follows from above remarks that the exponential map
$$\VEC t,k\mapsto e^{t_1\zeta}e^{t_2\zeta^3}\cdots e^{t_k\zeta^{2k-1}}
$$
gives an isomorphism between affine additive group of rank $\cc n+1,2$ and the
connected component of the identity $\Is_0$.

$\Is$ can be defined as an algebraic subvariety in the affine space $\aA$ by
equation
$$f^\vee\cdot f=f(-\zeta)f(\zeta)=1\,.
$$
If we use the coefficients $\OVEC a,n$ of polynomials
$$f(\zeta)=\BSER a,\zeta,k\in\aA
$$
as coordinates on $\aA$, then this equation is equivalent to the system of
$k+1$ quadratic equations
$$\left\{\begin{array}{rcl}
        a_0^2&=&1\\
       2a_2^2&=&-a_1^2\\
       2a_4^2&=&-2a_1a_3-a_2^2\\
              &\cdots&\\
      2a_{2k}&=&-2a_1 a_{2k-1}-2a_2a_{2k-2}-\cdots-2a_{k-1}a_{k+1}-a_k^2
    \end{array}\right.\quad.
$$
We see that for any choose of a value $a_0=\pm1$ and any fixed values of
$\cc{n+1},2$ odd coefficients $a_{2\nu+1}$ there exist a unique collection of
values of even coefficients such that $f(\zeta)$ is isometric. Hence, the
projection of the affine space $\aA$ onto the subspace generated by odd
coordinates gives an isomorphism of $\Is$ with the disjoint union of two such
subspaces.
\EP

Recall that in proof of the proposition (\ref{epscase}) we consider
non-degenerate bilinear form $(v,w)=\sp v,{\eta^nw}$ over the factor
$V/\ker(\eta^n)$. If the original bilinear form $\sp*,*$ is of type 1, then
this form $(*,*)$ is symmetric and factor $V/\ker(\eta^n)$ is $1$-dimensional.
Hence, the number $\varrho\bydef(v,v)=\sp v,{\eta^nv}$ modulo multiplication by
squers does not depend on $v\in V/\im(\eta)$.

\Pr\label{canone}
For any indecomposable rational form of type 1 on $(n+1)$-dimensional vector
space there exists a Jordan basis of $\zeta$ over quadratic extension
$\QQ(\sqrt{\varrho/{2^n}})$ such that the Gram matrix of the form at this basis
is equal to
$$\left(
  \begin{array}{cccccc}
        &&&&&1\\
        &0&&&-1&1\\
        &&&1&-1&\\
        &&-1&1&&\\
        &\Ddots&\Ddots&&0&\\
        (-1)^n&(-1)^{n-1}&&&&
  \end{array}
  \right)\quad.
$$
\EF
\proof
Let us fix an arbitrary Jordan basis
$$e_n\stackrel{\zeta}{\mapsto}
       e_{n-1}\stackrel{\zeta}{\mapsto}
        \cdots\stackrel{\zeta}{\mapsto}
           e_0\stackrel{\zeta}{\mapsto}e_{-1}=0
$$
for $\zeta$. We are going to find a selfdual operator
$$f=b_0+b_2\zeta^2+b_4\zeta^4+\cdots+b_{2k}\zeta^{2k}
$$
such that
$$\sp fe_i,{e_n}=\left\{\begin{array}{ccc}
                                       1&{\rm for}&i=0,1\\
                                       0&{\rm for}&i\ge2
                              \end{array}
                       \right.\;.
$$
In order to do this note that $\sp fv,{fw}=\sp v,{gw}$, where
$$g=f^*f=a_0+a_2\zeta^2+a_4\zeta^4+\cdots+a_{2k}\zeta^{2k}
$$
is selfdual too. Using orthogonality conditions (\ref{ortcond}) we obtain:
$$\sp e_{2\nu},{ge_n}=a_0\sp e_{2\nu},{e_n}+a_2\sp e_{2\nu},{e_{n-2}}+
                   a_4\sp e_{2\nu},{e_{n-4}}+\cdots+a_{2\nu}\sp
e_{2\nu},{e_{n-2\nu}}\;.
$$
Hence, we can choose the constants $\{a_2\nu\}$ such that
$$\sp e_{2\nu},{e_n}=\left\{\begin{array}{ccc}
                                       1&{\rm for}&\nu=0\\
                                       0&{\rm for}&\nu\ge1
                              \end{array}
                       \right.\;.
$$
Moreover, we can take $a_0=\sp e_o,{e_n}^{-1}$. In order to get $f$ from $g$ we
have to solve the system
$$\left\{\begin{array}{rcl}
        a_0&=&b_0^2\\
        a_2&=&2b_0b_2\\
        a_4&=&2b_0b_4+b_2^2\\
              &\cdots&\\
     a_{2k}&=&2b_0b_{2k}+2b_2b_{2k-2}+2b_4b_{2k-4}+\cdots
    \end{array}\right.\quad.
$$
It is possible over the qudratic extension $\QQ(\xi)$, where
$$\xi^2=\sp e_0,{e_n}=\varepsilon\sp e_n,{e_0}=\varepsilon\sp e_n,{\zeta^ne_n}=
                       \frac{1}{2^n}\sp e_n,{\eta^ne_n}=\frac{1}{2^n}\varrho\;.
$$

Finally, using formula (\ref{jogr}), for basic vectors with odd indices we get:
\begin{eqnarray*}
    2\sp e_{2\nu-1},{e_n}&=&\sp e_{2\nu-1},{e_n}-
                        \varepsilon\sp e_n,{e_{2\nu-1}}=
               \sp e_n,{\kappa e_{2\nu-1}}-\sp e_n,{\varepsilon e_{2\nu-1}}=\\
    &=&\sp e_n,{(\kappa-\varepsilon E) e_{2\nu-1}}=
                 2\varepsilon\bigl(\sp e_n,{e_{2\nu-2}}+
                                        \sp e_n,{e_{2\nu-3}}+\cdots\,\bigr)=\\
     &=&2\bigl(\sp {e_{2\nu-2}},{e_n}-\sp {e_{2\nu-3}},{e_n}+
                       \sp {e_{2\nu-4}},{e_n}-\cdots\,\bigr)\;,
\end{eqnarray*}
Hence,
$$\sp e_{2\nu-1},{e_n}=\left\{\begin{array}{ccc}
                                       1&{\rm for}&\nu=1\\
                                       0&{\rm for}&\nu\ge2
                              \end{array}
                       \right.
$$
and by (\ref{jogr}) we obtain the Gram matrix what is needed.
\EP
%
%
\No{{\normalsize $K_0(\PP_n)$} in more details}\label{Pn}
\SubNo{Notations}
In this section we consider in more details the module $\kK_n=K_0(\PP_n)$ with
the natural unimodular bilinear form
$$\sp E,F=\sum(-1)^\nu\dim\ext^\nu(E,F),
$$

Denote by $\pP_n\!\subset\!\QQ[t]$ the subspace of all polynomials of degree
$\le n$, and let $\mM_n\!\subset\!\pP_n$ be the \ZZ-submodule of all
polynomials taking integer values at all integer points. We will call such
polynomials {\it numerical\/}. Evidently, $\pP_n=\mM_n\otimes_\ZZ\QQ$.

The map
$$\MAP h,\kK_n,\mM_n\,\colon\;E\mapsto h_E(t)=\chi(E\otimes\oO(t))\;\mbox{(for
$t\bl\ZZ$)},
$$
which takes a coherent sheaf $E$ to its Hilbert polynomial $h_E(t)$, is an
isomorphism of \ZZ-modules. We identify $\kK_n$ with $\mM_n$ by this
isomorphism.

It is easy to check that this identification takes the \ZZ-basis
of $\kK_n$ consisting of the restrictions of the structure sheaf onto
subspaces:
$$\{\oO_{\PP_n},\oO_{\PP_{n-1}},\dots,\oO_{\PP_1},\oO_{\PP_0}\}
$$
to standard binomial \ZZ-basis
$$\{\gamma_n(t),\gamma_{n-1}(t),\dots,\gamma_0(t)\}
$$
of $\mM_n$ consisting of
\begin{eqnarray*}
       \gamma_{k}(t)&=&h_{\oO_{\PP_{k}}}(t)={t+k\choose
k}=\frac{1}{k!}(t+1)(t+2)\cdots(t+k)\\
                      &&\qquad\qquad\qquad\qquad\qquad\qquad\qquad\qquad
                                         {\rm for}\quad k=1,2,\dots,n\\
       \gamma_0(t)&=&h_{\oO_{\PP_0}}(t)\equiv1
\end{eqnarray*}

The restriction operator $E\mapsto E|_{\PP_{n-1}}$ onto a hyperplane
$\PP_{n-1}\!\subset\!\PP_n$ is represented in terms of $\mM_n$ by {\it left
difference operator\/}
$$\MAP{\nabla=1-e^{-D}},\mM_n,\mM_n\,\colon\;f(t)\mapsto\nabla f(t)\bydef
f(t)-f(t-1),
$$
where $D=d/dt$. Note that the polynomials $\gamma_\nu$ form a Jordan chain for
this operator, i.e.
$\nabla^m\gamma_\nu=\gamma_{\nu-m}$.

\SubNo{Canonical algebra}
We denote by $\aA_n$ the canonical algebra of all reflexive operators with
respect to the form on $\mM_n$ comming from $\kK_n$ under our identification.
Recall that it coincides with centralizer of $\kappa$ in
$\hom_\ZZ(\mM_n,\mM_n)$. To describe  $\aA_n$  we describe  first its
vectorisation $\aA_n\otimes\QQ$, i.e. the centralizer of $\kappa$ in
$\hom_\QQ(\pP_n,\pP_n)$.

The canonical operator of the bilinear form on $\kK_n$ coincides with the
Serre-Verdier dualizing operator, which takes a class of coherent sheaf $E$ to
a class $(-1)^nE(-n-1)$. Under the isomorphism $h$ this operator is identified
with the operator
$$\MAP{\kappa=(-1)^ne^{-(n+1)D}},\mM_n,\mM_n\,\colon\;f(t)\mapsto (-1)^n
f(t-n-1)
$$
Hence, the canonical operator has the form $\kappa=(-1)^n\Id+\eta$, where
\begin{equation}\label{eta}
\eta=(-1)^n(e^{-(n+1)D}-1)
\end{equation}
is nilpotent operator such that $\eta^n\ne0$, but $\eta^{n+1}=0$. So, in terms
of the previous paragraph, we get

\Pr $\mM_n\otimes\QQ$ is the space of type 1.
\EF\EP

In particular, the centralizer of $\kappa$ in $\hom_\QQ(\pP_n,\pP_n)$ is equal
to $\QQ[\eta]/\eta^{n+1}$. Since  $\QQ[\eta]/\eta^{n+1}=\QQ[D]/D^{n+1}$ by
(\ref{eta}), we get
\Cl $\aA_n\otimes\QQ=\QQ[D]/D^{n+1}$, where D=d/dt
\EF\EP

\Cl $\aA_n=\ZZ[\nabla]/\nabla^{n+1}$, where $\nabla=1-e^{-D}$.
\EF
\proof Of course, $\aA_n\otimes\QQ=\QQ[D]/D^{n+1}=\QQ[\nabla]/\nabla^{n+1}$.
So, we have to prove that an operator
$$\MAP{A=\BSER a,\nabla,n},\pP_n,\pP_n
$$
takes $\mM_n$ into $\mM_n$ if and only if all $a_\nu\bl\ZZ$. To do this we
apply $A$ to $\gamma_\nu$ and evaluate at the point $t=0$. We get
$a_\nu=A\gamma_\nu(0)$. Hence,
$$A\gamma_\nu\bl\mM_n\Leftrightarrow a_\nu\bl\ZZ.
$$
\EP

\SubNo{Tensoring and dualizing}
There are two more algebraic structures on the module $\kK_n$ --- the structure
of the ring with respect to the tensor product of locally free sheaves and the
involution $*$ taking a locally free sheaf $E$ to its dual $E^*={\cal H\it
om}(E,\oO)$. We carry these operations over the module $\mM_n$ by isomorphism
$h$ and denote by $\otimes$ and $*$ as well.

Since tensoring by the restriction of the structure sheaf onto hyperplane is
represented in terms of $\mM_n$ by the operator $\nabla$, we get immediately
that
$$\gamma_{n-\nu}\otimes\gamma_{n-\mu}=\gamma_{n-(\nu+\mu)}.
$$
Note, that $\gamma_n$ is the unit element with respect to tensor product.

Let us define the linear map $\MAP\ch^{-1},\aA_n,\mM_n$ by the rule
$$\ch^{-1}\,\colon\;A\mapsto A\gamma_n
$$
Evidently, this map is an isomorphism of \ZZ-modules. Moreover, the following
proposition holds:

\Pr $\ch^{-1}(AB)=\ch^{-1}(A)\otimes\ch^{-1}(B)$, i.e. the map $\ch^{-1}$
is an isomorphism of \ZZ-algebras, where the multiplication on $\mM_n$ is given
by the tensor product and the multiplication on $\aA_n$ is given by the
composition of operators (or multiplication of formal power series modulo
$\nabla^{n+1})$.
\EF
\proof It is sufficient to check the formula for basic operators $A=\nabla^k$,
$B=\nabla^m$, but in this case it is obvious.
\EP

The inverse isomorphism $\MAP\ch,\mM_n,\aA_n$ will be called the {\it Chern
character\/}. We identify $\mM_n$ (and $\kK_n$) with canonical algebra $\aA_n$
by this isomorphism and carry the involution $*$ over $\aA_n$ as well. We are
going to compare this involution with the involution $\scriptstyle\vee$, which
takes a reflexive operator to its dual with respect to the bilinear form in the
sense of (\ref{ca}).

\Pr For any $A=A(D)\bl\aA_n\otimes\QQ$ we have
$$A^\vee(D)=A^*(D)=A(-D).
$$
\EF
\proof The involution induced by the rule $D\mapsto-D$ takes the canonical
operator $\kappa$ to its dual $\kappa^{-1}$, and hence, this involution
coincides with the involution $\scriptstyle\vee$ of the canonical algebra.

Further, for {\it translation operator\/} $T=e^D\,\colon\;f(t)\mapsto f(t+1)$
we have $T^\vee=T^{-1}$. This means that $\scriptstyle\vee$ acts on $\mM_n$ by
taking
$$T^k\gamma_n(t)=\gamma_n(t+k)=h_{\oO_{\PP_n}(k)}(t)
$$
to
$$T^{-k}\gamma_n(t)=\gamma_n(t-k)=h_{\oO_{\PP_n}(-k)}(t).
$$
Hence, $\scriptstyle\vee$ coincides with $*$.
\EP

Since for any pair of locally free sheaves we have
$$\sp E,F=\chi(E^*\otimes F)=h_{E^*\otimes F}(0),
$$
we get for the scalar product on $\mM_n$ the formula
$$\sp f,g=f^*\otimes g\,(0).
$$
Hence, we obtain

\Cl In terms of operator $D$, the scalar product on $\aA_n$ is given by
$$\sp A(D),{B(D)}=A(-D)B(D)\gamma_n\,(0).
$$
\EF\EP

\SubNo{Standard basises and their Gram matrices}
We have seen in the previous paragraph that there exists a basis $\OSET\Xi,n$
over some quadratic extension $\aA_n\otimes\QQ(\xi)$
with Gram matrix
$$\left(
  \begin{array}{ccllll}
        &&&&&1\\
        &0&&&-1&1\\
        &&&1&-1&\\
        &&-1&1&&\\
        &\Ddots&\Ddots&&0&\\
        (-1)^n&(-1)^{n-1}&&&&
  \end{array}
  \right)\quad.
$$

Recall that the number $\xi$ was defined by quadratic equation
$$\xi^2=\frac{1}{2^n}\sp1,{\eta^n}=\left(\frac{n+1}{2}\right)^nD^n\gamma_n(0)
                                       =\left(\frac{n+1}{2}\right)^n\;.
$$
Hence, for even $n$ the basis in question is retional and for odd $n$ it exists
over
$\QQ(\sqrt{(n+1)/2})$.

Recall also that this basis coincides with some Jordan basis for the operator
$$\zeta=\frac{(-1)^n\kappa-E}{(-1)^n\kappa+E}=-\tanh\left(\frac{n+1}{2}D\right)
$$
and has a form
$\{\varphi\zeta^n,\,\varphi\zeta^{n-1},\,\dots\,,\,\varphi\zeta,\,\varphi\}$,
where $\varphi=\varphi(D)$ is an appropriate selfdual operator. Unfortunately,
I do not know
any general explicit formula for these basic operators. Indeed, it is not
difficult to calculate them for concrete $n$. For example, on
$K_0(\PP_2)\otimes\QQ$ we can take
$$\left\{\frac{3}{2}D^2\,,\;-D\,,\;\frac{2}{3}-\frac{1}{3}D^2\right\}
$$

We will suppose that some such basis is fixed in $\aA_n\otimes\QQ$ for each $n$
and will denote it by $\OSET\Xi,n$, where $\zeta\Xi_\nu=\Xi_{\nu-1}$.
Coordinates $\OSET z,n$ with respect to this basis may be considered as some
characteristic classes of sheaves and it would be very interesting to
investigate their geometrical sense.

An other rational basis of $\aA_n\otimes\QQ$, which is useful for calculations,
consists of
{\it Adams operators\/}
$$\Psi_k(D)\bydef\frac{\textstyle D^k}{\textstyle k!}\quad
                                \mbox{ (where $k=0,1,\dots n$).}
$$
Note, that for this basis we have $\Psi_k^*=(-1)^k\Psi_k$ too.

Let $A,B\bl\aA_n\otimes\QQ$ be decomposed by $\Psi_\nu$ as
$$A=\sum a_\nu\Psi_\nu\;,\quad B=\sum b_\nu\Psi_\nu.
$$
Then we can written
$$A^*B=\sum \alpha_\nu(A,B)\Psi_\nu,
$$
where each bilinear form $\alpha_k(A,B)$ depends on only first $(k+1)$
coefficients
$$\OVEC a,k\quad{\rm and}\quad\OVEC b,k
$$
of $A$, $B$. The precise expression for $\alpha_k$ is
$$\alpha_k(A,B)=\sum^k_{\nu=0}(-1)^\nu{k\choose\nu}a_\nu b_{k-\nu}.
$$
Note that $\alpha_k(A,B)$ is symmetric for even $k$ and skew-symmetric for odd
$k$.

\Pr\label{gramadams}
The original bilinear form on $\aA_n\otimes\QQ$ is decomposed by the forms
$\alpha_k(A,B)$ in the following way:
$$\sp A,B=\frac{1}{n!}\sum_{k=0}^n\sigma_{n-k}(1,2,\dots,n)\alpha_k(A,B),
$$
where
$\sigma_{n-k}(1,2,,\dots,n)=\sum\limits_{1\le\nu_1<\nu_2<\cdots<\nu_k\le
n}\nu_1\nu_2\cdots\nu_k$
is the value of $(n-k)$-th elementary symmetrical polynomial at the integer
point $(1,2,\dots,n)$.
\EF
\proof
It is easy to check that
$$\frac{D^k}{k!}(t+1)(t+2)\cdots(t+n)=\sigma_{n-k}((t+1),(t+2),,\dots,(t+n)).
$$
Hence,
\begin{eqnarray*}
\sp
A,B&=&A^*B\gamma_n\,(0)=\sum\alpha_k(A,B)\frac{D^k}{k!}\frac{1}{n!}(t+1)(t+2)\cdots(t+n)|_{t=0}\\
       &=&\frac{1}{n!}\sum\sigma_{n-k}(1,2,,\dots,n)\alpha_k(A,B).
\end{eqnarray*}
\EP

For example, Gram matrices of Adams basises for $\PP_2$, $\PP_3$, $\PP_4$ are
the following:
$$\frac{1}{2}
 \left(
    \begin{array}{rrr}
         2&3&1\\
         -3&-2&0\\
         1&0&0
     \end{array}
  \right)\;,\qquad\qquad\qquad
 \frac{1}{6}
\left(
    \begin{array}{rrrr}
         6&11&6&1\\
         -11&-12&-3&0\\
         6&3&0&0\\
         -1&0&0&0
     \end{array}
  \right)\;,
$$
$$
 \frac{1}{24}
\left(
    \begin{array}{rrrrr}
         24&50&35&10&1\\
         -50&-70&-30&-4&0\\
         35&30&6&0&0\\
         -10&-4&0&0&0\\
         1&0&0&0&0
     \end{array}
  \right)\;.
$$

\SubNo{Isometries and their invariants}
We have seen in the previous paragraph that the isometry group $\Is$ has two
connected components and the component of the identity coincides with the image
of the exponential map applyed to the subspace $\aA_-\otimes\QQ$ of all
antiselfdual operators. In terms of operator $D$, this subspace consists of all
operators $A(D)$ represented by odd power series.

\Pr Two nondegenerate operators $A,B\in\aA_n\otimes\QQ$ belong to the same
orbit of the natural action of $\Is$ on $\aA_n\otimes\QQ$ by multiplications if
and only if $A^*A=B^*B$ in  $\aA_n\otimes\QQ$.
\EF
\proof
If $A=B\Phi$ for some $\Phi\bl\Is$, then
$$A^*A=B^*B\Phi^*\Phi=B^*B,
$$
because $\Phi^*\Phi=1$.

At the same time, if $A^*A=B^*B$ and $A$ and $B$ are invertible, then
$\Phi=AB^{-1}$ satisfies the condition
$$\Phi^*=A^*{B^*}^{-1}=BA^{-1}=\Phi^{-1},
$$
and hence, it is isometric.
\EP

\SubNo{The rank}
Let us define the {\it rank functional\/}
$$\MAP\rk,\aA_n,\ZZ
$$
by the rule:
$$\rk(A)\bydef\sp
A,{\nabla^n}=\varepsilon\sp\nabla^n,A\;\fain{A=A(\nabla)},{\aA_n}\;.
$$
It follows from the orthogonality conditions (\ref{ortcond}) that
$$\rk(\BSER x,\nabla,n)=x_0\;.
$$
Hence, $\forall\,A,B$ we have $\rk(A)\rk(B)=\sp
A,{\nabla^nB}=\varepsilon\sp\nabla^nA,B$.

Geometrically, if the operator $A$ corresponds to the class of a locally free
sheaf, then $\nabla^nA$ corresponds to its restriction onto a point. Hence, the
rank defined above coincides in this case with the usual rank of locally free
sheaf.

Since $\nabla=1-e^{-D}$ and $\nabla^n=\varepsilon D^n$ we can calculate rank in
terms of $D$ by the formula
$$(\rk A)^2=\varepsilon\sp A,{D^nA}\;.
$$
So, $\rk(a_0\Psi_0+a_1\Psi_1+\cdots+a_n\Psi_n)=a_0$.

In terms of an other useful operator $\zeta=-\tanh((n+1)D/2)$ we have
$$D=-\frac{1}{n+1}\log\left(\frac{1+\zeta}{1-\zeta}\right)\quad\mbox{ and
}\quad
                        D^n=\left(\frac{-2}{n+1}\right)^n\zeta^n\;.
$$
Hence,
$$(\rk A)^2=\left(\frac{2}{n+1}\right)^n\varepsilon\sp A,{\zeta^n A}\;.
$$
So,
$\Bigl(\rk(z_0\Xi_0+z_1\Xi_1+\cdots+z_n\Xi_n)\Bigr)^2=\left(\frac{\textstyle
2}{\textstyle n+1}\right)^nz _n^2$.

\SubNo{$K_0(\PP_2)$ and Markov chain}
Let $\mM$ be a free \ZZ-module of rank 3 equipped with an unimodular integer
bilinear form. Jordan normal form of corresponding the canonical operator on
$\mM\otimes\CC$ may be only one of the following:
$$\left(
    \begin{array}{ccc}
         1&1&0\\
         0&1&1\\
         0&0&1
     \end{array}
  \right)\;,\qquad
  \left(
    \begin{array}{rrc}
         -1&1&0\\
         0&-1&0\\
         0&0&1
     \end{array}
  \right)\;,\qquad
 \left(
    \begin{array}{ccc}
         \lambda&0&0\\
         0&\lambda^{-1}&0\\
         0&0&1
     \end{array}
  \right)\;,
$$
where $\lambda\ne\pm1$. These three cases are distinguished by a value of the
trace $\tr(\kappa)$, which is equal to 3, $-1$ and $1+\lambda+\lambda^{-1}$,
where $\lambda\ne\pm1$, respectively.

Suppose now that the form on $\mM$ admits some semiorthonormal basis with Gram
matrix
\begin{equation}\label{abc}
 \chi=
  \left(
    \begin{array}{ccc}
         1&a&b\\
         0&1&c\\
         0&0&1
     \end{array}
  \right)\;.
\end{equation}
Easy computation gives $\tr(\kappa)=\tr(\chi^{-1}\chi^t)=3-a^2-b^2-c^2+abc$. We
get

\Pr An integer bilinear form (\ref{abc}) is of type 1 (in the sense of previous
paragraph) if and only if the numbers $\{a,b,c\}$ satisfy the {\it tripled
Markov equation\/}:
$$a^2+b^2+c^2=abc\;.
$$
\EF\EP

It is well known (see [Ca]) that all solutions of tripled Markov equations are
obtained from
the initial solution $\{3,3,3\}$ by use of the following two procedures:
\Cs
\cond changing signs of any two numbers;
\cond changing a value of one of numbers via Vieta theorem:
$$a\mapsto bc-a\,,\mbox{ or }b\mapsto ac-b\,,\mbox{ or }c\mapsto ab-c\,.
$$
\EC

Note now that we can change signs of any two elements of Gram matrix by
changing a sign of one of basic vectors. Further, the following three mutations
of a semiorthonormal basis $\{e_0,e_1,e_2\}$ (see \ref{mutations}):
\begin{eqnarray*}
   L_1\,\colon\;\{e_0,e_1,e_2\}&\longmapsto&\{e_1-\sp
e_0,{e_1}e_0\,,\;e_0\,,\;e_2\}\\
   L_2\,\colon\;\{e_0,e_1,e_2\}&\longmapsto&\{e_0\,,\;e_2-\sp
e_1,{e_2}e_1\,,\;e_1\}\\
   R_2\,\colon\;\{e_0,e_1,e_2\}&\longmapsto&\{e_0\,,\;e_2\,,\;e_1-\sp
e_1,{e_2}e_2\}
\end{eqnarray*}
change Gram matrix by the rules
\begin{eqnarray*}
  \left(
    \begin{array}{ccc}
         1&a&b\\
         0&1&c\\
         0&0&1
     \end{array}
  \right)
 &\longmapsto&
  \left(
    \begin{array}{crc}
         1&-a&c-ab\\
         0&1&b\\
         0&0&1
     \end{array}
  \right)\\
  \left(
    \begin{array}{ccc}
         1&a&b\\
         0&1&c\\
         0&0&1
     \end{array}
  \right)
 &\longmapsto&
  \left(
    \begin{array}{ccr}
         1&b-ac&a\\
         0&1&-c\\
         0&0&1
     \end{array}
  \right)\\
\left(
    \begin{array}{ccc}
         1&a&b\\
         0&1&c\\
         0&0&1
     \end{array}
  \right)
 &\longmapsto&
  \left(
    \begin{array}{ccc}
         1&b&a-bc\\
         0&1&-c\\
         0&0&1
     \end{array}
  \right)\;.
\end{eqnarray*}

Hence, we have proved

\Pr
Any semiorthonormal basis of  an integer bilinear form of type 1 can be
transformed using the braid group action and changing signs of basic vectors to
the semiorthonormal basis with Gram matrix
$$\left(
    \begin{array}{ccc}
         1&3&3\\
         0&1&3\\
         0&0&1
     \end{array}
  \right)\;.
$$
\EF\EP

The last matrix coincides with the Gram matrix of the basis $\{\oO, {\cal
T}(-1),\oO(1)\}$ of $K_0(\PP_2)$, where ${\cal T}(-1)$ is the twisted tangent
sheaf. We get

\Cl There exists a unique up to integer isometries integer bilinear form of
type 1, which admits a semiorthonormal basis. This form coincides with the form
on $K_0(\PP_2)$.
\EF\EP
%
%
\vspace{3ex}\par\noindent
\begin{flushleft}
 {\large\bf References.}
\end{flushleft}
\vspace{1ex}
\begin{description}
\item{[Bo1]:}
A.I.Bondal.
{\it Representations of associative algebras and coherent sheaves.\/}
Math. USSR Izv, 34 (1990) no1, p.23-42 (english translation).
\item{[Bo2]:}
A.I.Bondal.
{\it Helices, representations of quivers and Koszul algebras.\/}
In: {\it Helices and Vector Bundles\/}, p.75-95, Lond.Math.Soc. L.N.S. 148
(1990).
\item{[BoKa]:}
A.I.Bondal, M.M.Kapranov
{\it Representable functors, Serre functors and mutations.\/}
Math. USSR Izv, 35 (1990) no3, p.519-541 (english translation).
\item{[Ca]:}
J.W.S.Cassels.
{\it An introduction to Diophantine approximation.\/}
Hafner Publishing Company, New York, 1972.
\item{[Go1]:}
A.L.Gorodentsev.
{\it Surgeries of exceptional vector bundles on $\PP_n$.\/}
Math. USSR Izv,  32 (1989) no1, p.1-13 (english translation).
\item{[Go2]:}
A.L.Gorodentsev.
{\it Exceptional vector bundles on surfaces with moving anticanonical
divisor.\/}
Math. USSR Izv, 33 (1989) no1, p.67-83 (english translation).
\item{[Go3]:}
A.L.Gorodentsev.
{\it Exceptional objects and mutations in derived categories.\/}
in: {\it Helices and Vector Bundles\/}, p.57-75, Lond.Math.Soc. L.N.S. 148
(1990).
\item{[Go4]:}
A.L.Gorodentsev.
{\it Helix theory and non-symmetric bilinear forms.\/}
In: {\it Algebraic Geometry and its Applications\/} (proceedings of 8th
Alg.Geom.Conf., Yaroslavl'1992). Aspects of Math. (1994).
\item{[GoRu]:}
A.L.Gorodentsev, A.N.Rudakov.
{\it Exceptional vector bundles on projective spaces.\/}
Duke Math.J. 54 (1987), p.115-130.
\item{[HoPe]:}
W.V.D.Hodge, D.Pedoe.
{\it Methods of Algebraic Geometry. Vol.1.\/}
Cambridge Univ. Press, 1953.
\item{[Ma]:}
A.I.Mal'zev.
{\it Foundations of Linear Algebra\/} (in russian).
OGIZ, Moskva-Leningrad, 1948.
\item{[No]:}
D.Yu.Nogin.
{\it Helices of period 4 and Markov-type equations.\/}
Math. USSR Izv, 37 (1991) no1, p.209-226 (english translation).
\item{[Ru]:}
A.N.Rudakov
{\it Integer-valued bilinear forms and vector bundles.\/}
Math. of the USSR Sbornik, 66 (1990) no1, p.189-197 (english translation).
\item{[Se]:}
J.-P.Serre.
{\it Lie algebras and Lie Groups.}
, Benjamin, 1965.
\end{description}
\end{document}